\documentclass[a4paper,twoside,10pt,reqno]{article}
\usepackage[english]{babel}
\usepackage[dvips]{graphicx}
\usepackage[T1]{fontenc}
\usepackage[latin1]{inputenc}
\usepackage{amssymb,amsthm,dsfont,mathrsfs}
\usepackage{amsfonts,amsmath,euscript}
\usepackage{color}

\newtheorem{defi}{Definition}[section]
\newtheorem{teo}{Theorem}[section]
\newtheorem{pro}[teo]{Proposition}
\newtheorem{cor}[teo]{Corollary}

\newtheorem{rem}{Remark}[section]
\newtheorem{ex}{Example}[section]

\newcounter{example}[section]

           \addto{\captionsenglish}{}

\newcommand{\er}{\mathbb{R}}

\newcommand{\hs}{\hspace{3pt}}
\newcommand{\dst}{\displaystyle}
\newcommand{\dem}{{\bf Dem. }}
\newcommand{\fdem}{$\square$}

\newcommand{\nid}{\noindent}
\newcommand{\nn}{\nonumber}

\newcommand{\titulo}[1]{\mbox{} \\ \noindent \textit{\textbf{\Large #1}}\\}

\renewcommand{\abstract}[1]{{\small \noindent \textbf{Abstract:} #1\\}}

\newcommand{\pchave}[1]{{\small \noindent \textbf{Keywords:} #1\\}}

\begin{document}

\titulo{Fragility Index of block tailed vectors}
%
\vspace{0.5cm}

\textbf{Helena Ferreira} Department of Mathematics, University of
Beira
Interior, Covilhã, Portugal\\

\textbf{Marta Ferreira} Department of Mathematics, University of
Minho, Braga, Portugal\\

\abstract{Financial crises are a recurrent phenomenon with important
effects on the real economy. The financial system is inherently
fragile and it is therefore of great importance to be able to
measure and characterize its systemic stability. Multivariate
extreme value theory provide us such a framework through the
\emph{fragility index} (Geluk \cite{gel+}, \emph{et al.}, 2007; Falk
and Tichy, \cite{falk+tichy1,falk+tichy2} 2010, 2011). Here we
generalize this concept and contribute to the modeling of the
stability of a stochastic system divided into blocks. We will find
several relations with well-known tail dependence measures in
literature, which will provide us immediate estimators. We end with
an application to financial data.}

\pchave{multivariate extreme value theory, tail dependence,
fragility index, extremal coefficients}

\section{Introduction}


{\color{blue}In the last decade, dependencies between financial
asset returns have increased, mostly as a consequence of
globalization effects and relaxed market regulation. Therefore, the
concept of tail dependence has been discussed in financial
applications related to market or credit risk, e.g., Hauksson
\emph{et al.} \cite{hauk+} (2001), Ané and Kharoubi \cite{ane}
(2003), Junker and May \cite{junk} (2005), Embrechts and Puccetti
\cite{emb+puc} (2010).
%
%
%
The natural framework to model extremal dependence turns out to be
the multivariate extreme value theory. The study of systemic
stability is an important issue within this context of extreme risk
dependence. The fragility of a system have been addressed to the
\emph{Fragility Index} (FI) introduced in Geluk \emph{et
al.}\hs(\cite{gel+}, 2006).} More precisely, consider a random
vector $\mathbf{X}=(X_1,...,X_d)$ and
$N_x:=\sum_{i=1}^d\mathds{1}_{\{X_i>x\}}$ the number of exceedances
among $X_1,...,X_d$ above a threshold $x$. The FI corresponding to
$\mathbf{X}$ is the asymptotic conditional expected number of
exceedances, given that there is at least one exceedance, i.e., $FI
= \lim_{x\to\infty} E(N_x | N_x > 0)$. The stochastic system
$\{X_1,...,X_d\}$ is called fragile whenever $FI>1$. Theoretical
developments, namely, the asymptotic distribution of $N_x$
conditional to $N_x>0$ can be seen in Falk and Tichy
(\cite{falk+tichy1,falk+tichy2}, 2010, 2011).

In this work we generalize some properties of the FI presented in
the references above, contributing to the modeling of the
stability of a stochastic system divided into blocks.\\

We shall state some notation that will be used throughout the paper.\\
Consider $\mathcal{D}=\{I_1,...,I_s\}$ a partition of
$D=\{1,...,d\}$. For the random vector $\mathbf{X}=(X_1,...,X_d)$,
let $\mathbf{X}_{I_j}$ be a sub-vector of $\mathbf{X}$ whose
components have indexes in $I_j$, with $j=1,...,s$. If $F$ denotes
the d.f.\hs of $\mathbf{X}$ then $F_{I_j}$ denotes the d.f.\hs of
sub-vector $\mathbf{X}_{I_j}$, $j=1,...,s$, and $F_i$ the marginal
d.f., $i=1,...,d$. Let $\mathbf{x}_{I_j}$ be a vector of length
$|I_j|$ with components equal to $x\in\er$. We will say that
$\mathbf{X}_{I_j}$ is the $j^{th}$ block of random vector
$\mathbf{X}$ and denote by $N_{\mathbf{x}}$ the number of blocks
where it occurs at least one exceedance of $\mathbf{x}$, i.e.,
$$
N_\mathbf{x}=\sum_{j=1}^s\mathds{1}_{\{\mathbf{X}_{I_j}\not\leq \mathbf{x}_{I_j}\}}.
$$
All operations and inequalities on vectors are meant componentwise.

\begin{defi}
The Fragility Index (FI) of a random vector
$\mathbf{X}=(X_1,...,X_d)$ relative to partition $\mathcal{D}$ is
\begin{eqnarray}\label{fi}
\lim_{x\to\infty}E(N_{\mathbf{x}}|N_{\mathbf{x}}>0),
\end{eqnarray}
whenever the limit exists and is denoted
$FI(\mathbf{X},\mathcal{D})$.
\end{defi}

If we consider $I_j=\{j\}$,  $j=1,...,d$, we find the FI introduced
in Geluk \emph{et al.} (\cite{gel+}, 2007) and latter study by Falk
and Tichy (\cite{falk+tichy1,falk+tichy2}, 2010, 2011). This
partition will be denoted as $\mathcal{D}^*$, i.e.,
$\mathcal{D}^*=\{I_j=\{j\}: j=1,...,d\}$. \\

 We give particular emphasis to random vectors in the domain
of attraction of a multivariate extreme value distribution (MEV) and
consider either the case of identically distributed margins or  tail
equivalent margins in the sense considered in Falk and Tichy
(\cite{falk+tichy2}, 2011). In Section \ref{sapNx} we present some
asymptotic properties of the distribution of $N_x$ conditional to
$N_x>0$, and find generalizations of results in Falk and Tichy
(\cite{falk+tichy2}, 2011). We prove that
$FI(\mathbf{X},\mathcal{D})$ exists and relates with the extremal
coefficients $\epsilon$ of Tiago de Oliveira (\cite{tiago}, 1962/63)
and Smith (\cite{smith}, 1990) in case of identically distributed
margins (Section \ref{sFIblocs}). We define generalized versions of
the multivariate tail dependence coefficients of Li (\cite{li3},
2009) and extend some of its results (Section \ref{sli}). In Section
\ref{sli+fi} we relate these latter coefficients with
$FI(\mathbf{X},\mathcal{D})$.

For independent margins we have an unit FI. However, the stability
of a stochastic system at higher levels can also be characterized by
asymptotic independence (Geluk \emph{et al.}, \cite{gel+}, 2007).
Asymptotic independence means that the dependency when present
vanishes at extreme quantiles and the system is said to be weakly
fragile, albeit possibly correlated (e.g., gaussian vectors). We
extend the concept of asymptotic independent FI  in Geluk \emph{et
al.}\hs(\cite{gel+}, 2007) for blocks. A second measure is also
presented by extending the $2$-blocks asymptotic independent
coefficient in Ferreira and Ferreira (\cite{hf+mf3}, 2011) to the
case of $s$-blocks, with $s>2$. This issue is considered in Section
\ref{saifi}.

Our results relating $FI(\mathbf{X},\mathcal{D})$ with well-known
tail dependence measures, for which estimators and respective
properties have already been study in literature, will provide us
immediate estimators (Section \ref{sestim}). We end with an
application to financial data.


\section{Asymptotic Properties of $N_x$}\label{sapNx}

In this section we present some asymptotic properties of the
distribution of $N_x$ conditional to $N_x>0$. We start to relate
this latter with
 $\epsilon_A^G$, the extremal coefficients (Tiago de
Oliveira, \cite{tiago} 1962/63; Smith, \cite{smith} 1990) of the
sub-distribution functions of the MEV $G$ corresponding to margins
in $A$, i.e., by assuming $G$ has unit Fréchet margins,
\begin{eqnarray}\label{extremalcoef}
\epsilon_{A}^G=-\log
G(\mathds{1}_A^{-1}(1),...,\mathds{1}_A^{-1}(d))
\end{eqnarray}
 where, for all
$A\subset D$, $x\in\er$,
$$
\mathds{1}_A^{-1}(x)=\left\{
\begin{array}{ll}
           1&,\,x\in A\\
           \infty &,\, x\not\in A.
\end{array} \right.
$$
 They may be written through the stable tail dependence
function (Huang, \cite{huang} 1992):
\begin{eqnarray}\label{stdf}
l_G({x_1^{-1},...,x_d^{-1}})=-\log G(x_1,...,x_d)=-\log C_G(e^{-1/x_1},...,e^{-1/x_d}),
\end{eqnarray}
where {} $C_G$ is the copula of $G$, i.e.,
\begin{eqnarray}\label{copula}
G(x_1,...,x_d)=C_G(G_1(x_1),...,G_d(x_d)), \,\,\,(x_1,...,x_d)\in \mathbb{R}^d.
\end{eqnarray}

In the sequel we will use notation $\mathcal{I}(A)=\cup_{j\in
A}I_j$.

\begin{pro}\label{p1}
If $\mathbf{X}$ has d.f.\hs $F$ with identically distributed
continuous margins and belongs to the domain of attraction of a MEV
$G$ with unit Fréchet margins then, for each $k\in\{1,...,s\}$, we
have
$$
\dst\lim_{x\to\infty}P(N_x=k|N_x>0)=\frac{1}{\epsilon_D^G}\sum_{S\subset\{1,...,s\};|S|=k}
\sum_{T\subset S}(-1)^{|T|+1}\epsilon_{\mathcal{I}(T\cup S^C)}^G
$$
\end{pro}
\dem We have, successively,
\begin{eqnarray}\label{p1.1}
\begin{array}{rl}
P(N_x=k|N_x>0)= &
\frac{1}{1-P(N_x=0)}\sum_{S\subset \{1,...,s\};|S|=k}P(\cap_{j\in S}\mathbf{X}_{I_j}\not\leq \mathbf{x}_{I_j},
\cap_{j\not\in S}\mathbf{X}_{I_j}\leq \mathbf{x}_{I_j})\vspace{0.35cm}\\
=&\frac{1}{1-P(N_x=0)}\sum_{S\subset \{1,...,s\};|S|=k}\sum_{T\subset S}(-1)^{|T|}P(
\cap_{j\in T\cup S^C}\mathbf{X}_{I_j}\leq \mathbf{x}_{I_j})\vspace{0.35cm}\\
=&\frac{1}{1-F(\mathbf{x})}\sum_{S\subset \{1,...,s\};|S|=k}\sum_{T\subset S}(-1)^{|T|+1}
(1-F(x\mathds{1}^{-1}_{\mathcal{I}(T\cup S^c)}(1),...,x\mathds{1}^{-1}_{\mathcal{I}(T\cup S^c)}(d)))
\end{array}
\end{eqnarray}
since $\sum_{T\subset S}(-1)^{|T|+1}=0$. Assuming w.l.o.g.\hs that
$F$ has unit Pareto marginals, we obtained
\begin{eqnarray}\label{p1.2}
\begin{array}{rl}
&P(N_x=k|N_x>0)=  \vspace{0.35cm}\\
=&\frac{1}{1-C_F(\mathbf{1}-\frac{1}{x}\mathbf{1})}\sum_{S\subset \{1,...,s\};|S|=k}\sum_{T\subset S}(-1)^{|T|+1}
(1-C_F(1-\frac{1}{x}\mathds{1}^{-1}_{\mathcal{I}(T\cup S^c)}(1),...,1-\frac{1}{x}\mathds{1}^{-1}_{\mathcal{I}(T\cup S^c)}(d))).
\end{array}
\end{eqnarray}

By hypothesis, $F$ belongs to the domain of attraction of a MEV $G$,
which is equivalent to (de Haan and de Ronde, \cite{haan+ronde}
1998):
\begin{eqnarray}\label{p1.3}
\lim_{t\to\infty}\frac{1-C_F(1-y_1/x,...,1-y_d/x)}{1/x}=-\log C_G(e^{-y_1},...,e^{-y_d}),\,(y_1,...,y_d)\geq \mathbf{0}.
\end{eqnarray}

Taking limits in (\ref{p1.1}) and dividing both members by $1/x$,
conditions (\ref{p1.2}) and (\ref{p1.3}) lead us to
\begin{eqnarray}\nn
\begin{array}{rl}
&\lim_{x\to\infty}P(N_x=k|N_x>0)= \vspace{0.35cm}\\
=&\frac{1}{-\log C_G(e^{-1},...,e^{-1})}\sum_{S\subset \{1,...,s\};|S|=k}\sum_{T\subset S}(-1)^{|T|+1}(-\log
C_G(e^{-\mathds{1}_{\mathcal{I}(T\cup S^c)}(1)},...,e^{-\mathds{1}_{\mathcal{I}(T\cup S^c)}(d)}))\vspace{0.35cm}\\
=&\frac{1}{\epsilon_D^G}\sum_{S\subset \{1,...,s\};|S|=k}\sum_{T\subset S}(-1)^{|T|+1}\epsilon_{\mathcal{I}(T\cup S^C)}^G.\,\, $\fdem$
\end{array}
\end{eqnarray}

The previous result can be generalized to random vectors
$\mathbf{X}$ with equivalent marginal distributions, in the sense
that, there exists a d.f.\hs $H$ such that,
\begin{eqnarray}\label{equivmar}
\lim_{x\to w(H)}\frac{1-F_i(x)}{1-H(x)}=\gamma_i\in(0,\infty),\,i=1,...,d,
\end{eqnarray}
where $w(H)$ is the right-end-point of $H$. In this case it is no
longer possible an interpretation based on extremal coefficients, as
can be seen in the following result.

\begin{pro}\label{p2}
If $\mathbf{X}$ has d.f.\hs $F$ with equivalent marginal
distributions in the sense of (\ref{equivmar}), and belongs to the
domain of attraction of a MEV $G$ with unit Fréchet margins then,
for each $k\in\{1,...,s\}$, we have
\begin{eqnarray}\nn
\begin{array}{rl}
&\lim_{x\to w(H)}P(N_x=k|N_x>0)= \vspace{0.35cm}\\
=&\frac{1}{\log C_G(e^{-\gamma_1},...,e^{-\gamma_d})}\sum_{S\subset \{1,...,s\};|S|=k}\sum_{T\subset S}(-1)^{|T|+1}
\log C_G(e^{-\gamma_1\mathds{1}_{\mathcal{I}(T\cup S^c)}(1)},...,e^{-\gamma_d\mathds{1}_{\mathcal{I}(T\cup S^c)}(d)})
\end{array}
\end{eqnarray}
\end{pro}
\dem Observe that
\begin{eqnarray}\label{p2.1}
\begin{array}{rl}
&1-F(x\mathds{1}^{-1}_{\mathcal{I}(T\cup S^c)}(1),...,x\mathds{1}^{-1}_{\mathcal{I}(T\cup S^c)}(d))\vspace{0.35cm}\\
=& 1-C_F(1-(1-F_1(x\mathds{1}^{-1}_{\mathcal{I}(T\cup S^c)}(1))),...,1-(1-F_1(x\mathds{1}^{-1}_{\mathcal{I}(T\cup S^c)}(d))))\vspace{0.35cm}\\
=& 1-C_F(1-t_1^{\mathcal{I}(T\cup S^c)}(x),...,  1-t_d^{\mathcal{I}(T\cup S^c)}(x))
\end{array}
\end{eqnarray}
where
$$
t_i^A(x)=\mathds{1}_A(i)\frac{1-F_i(x)}{1-H(x)}(1-H(x)),\,\,i=1,...d\,.
$$
{Applying (\ref{equivmar}), we have
\begin{eqnarray}\label{p2.2}
\begin{array}{rl}
&\lim_{x\to w(H)}\frac{1-C_F(1-t_1^{\mathcal{I}(T\cup S^c)}(x),...,  1-t_d^{\mathcal{I}(T\cup S^c)}(x))}{1-H(x)} \vspace{0.35cm}\\
=&-\log C_G(e^{-\gamma_1\mathds{1}_{\mathcal{I}(T\cup S^c)}(1)},...,e^{-\gamma_d\mathds{1}_{\mathcal{I}(T\cup S^c)}(d)}).
\end{array}
\end{eqnarray}
The result follows by retaking expression in (\ref{p1.1}) and
considering (\ref{p2.1}) and (\ref{p2.2}). \fdem}\\

If in particular we consider $I_j=\{j\}$, $j=1,...,d$, we find the
result of Falk and Tichy (\cite{falk+tichy2}, 2011). \\

The results above can also be obtained through the relation between
$C_G$ and $D$-norms presented in Aulbach \emph{et al.} (\cite{au},
2011). However, we have chosen to present self-contained proofs
using the usual arguments of multivariate extreme value theory that
are more familiar.\\

\section{The Fragility Index for blocks}\label{sFIblocs}

In this section we compute the FI for blocks given in (\ref{fi}),
whenever $\mathbf{X}$ has equally distributed or tail equivalent
margins in the sense of (\ref{equivmar}), belonging to the domain of
attraction of a MEV $G$ with unit Fréchet margins.

\begin{pro}\label{p3}
If $\mathbf{X}$ has d.f.\hs $F$ with identically distributed
continuous margins and belongs to the domain of attraction of a MEV
$G$ with unit Fréchet margins, we have
$$
FI(\mathbf{X},\mathcal{D})=\frac{\sum_{j=1}^s
\epsilon_{I_j}^G}{\epsilon_D^G}\,.
$$
\end{pro}
\dem Observe that
\begin{eqnarray}\nn
\begin{array}{rl}
FI(\mathbf{X},\mathcal{D})= &\sum_{j=1}^s\lim_{x\to\infty}P(\mathbf{X}_{I_j}\not \leq \mathbf{x}_{I_j}|N_x>0) \vspace{0.35cm}\\
=&\sum_{j=1}^s\lim_{x\to\infty}\frac{1-F(x\mathds{1}^{-1}_{I_j}(1),...,x\mathds{1}^{-1}_{I_j}(d))}{1-F(x,...,x)}
\vspace{0.35cm}\\
=&\sum_{j=1}^s\lim_{x\to\infty}\frac{1-C_F(1-\frac{1}{x}\mathds{1}_{I_j}(1),...,1-\frac{1}{x}\mathds{1}_{I_j}(d))/(1/x)}
{1-C_F(1-\frac{1}{x},...,1-\frac{1}{x})/(1/x)}\vspace{0.35cm}\\
=& \frac{1}{\epsilon_D^G}\sum_{j=1}^s
\epsilon_{I_j}^G\,.\,\,  $\fdem$
\end{array}
\end{eqnarray}
In particular, for partitions corresponding to the family of all
margins we obtain the known result
$FI(\mathbf{X},\mathcal{D}^*)=d/\epsilon_D^G$. Hence, the FI of the
system divided into blocks is smaller than the system itself, i.e.,
$FI(\mathbf{X},\mathcal{D})\leq FI(\mathbf{X},\mathcal{D}^*)$.

\begin{rem}
We can relate $FI(\mathbf{X},\mathcal{D})$ with the fragility
indexes of the whole system and of each block. More precisely,
$FI(\mathbf{X},\mathcal{D})$ is a convex linear combination of the
ratios
$FI(\mathbf{X},\mathcal{D}^*)/FI(\mathbf{X}_{I_j},\mathcal{D}^*_{I_j})$,
since we can write
$$
FI(\mathbf{X},\mathcal{D})=\sum_{j=1}^s\frac{|I_j|}{d}\frac{FI(\mathbf{X},\mathcal{D}^*)}{FI(\mathbf{X}_{I_j},\mathcal{D}^*)}.
$$
Furthermore, it is also a weighted mean of those ratios:
$$
FI(\mathbf{X},\mathcal{D})=\frac{1}{s}\sum_{j=1}^s\frac{|I_j|s}{d}\frac{FI(\mathbf{X},\mathcal{D}^*)}{FI(\mathbf{X}_{I_j},\mathcal{D}^*)}.
$$
\end{rem}

\begin{pro}[inter-blocks dependence]\label{p4}
Under the conditions of Proposition \ref{p3}, 
we have
\begin{itemize}
\item[(i)] $1\leq FI(\mathbf{X},\mathcal{D})\leq \frac{\sum_{j=1}^s
\epsilon_{I_j}^G}{\bigvee_{j=1}^s\epsilon_{I_j}^G}.$
\item[(ii)] $FI(\mathbf{X},\mathcal{D})=1$ if and only if
$\mathbf{X}_{I_j}$, $j=1,...s$ are independent random vectors.
\item[(iii)] $FI(\mathbf{X},\mathcal{D})=\frac{\sum_{j=1}^s
\epsilon_{I_j}^G}{\bigvee_{j=1}^s\epsilon_{I_j}^G}$ if and only
if $\mathbf{X}_{I_j}$, $j=1,...s$ are totally dependent random
vectors in the sense
$G(\mathbf{x})=\bigwedge_{j=1}^sG_{I_j}(\mathbf{x}_{I_j})$.
\end{itemize}
\end{pro}
\dem By the inequalities
$$
\prod_{j=1}^sG_{I_j}(\mathbf{x}_{I_j})\leq G(\mathbf{x})\leq \bigwedge_{j=1}^s G_{I_j}(\mathbf{x}_{I_j})
$$
we have, successively,
$$
\begin{array}{c}
\dst (G_{\{1\}}(x))^{ \sum_{j=1}^s\epsilon_{I_j}^G}\leq G_{\{1\}}^{\epsilon_{D}^G}(x)\leq  (G_{\{1\}}(x))^{\bigvee_{j=1}^s\epsilon_{I_j}^G}\vspace{0.35cm}\\
\dst\bigvee_{j=1}^s\epsilon_{I_j}^G\leq \epsilon_{D}^G\leq  \sum_{j=1}^s\epsilon_{I_j}^G\,.\,\,$\fdem$
\end{array}
$$

\begin{pro}[intra-blocks dependence]\label{p5}
Under the conditions of Proposition \ref{p3}, 
we have
\begin{itemize}
\item[(i)] $\frac{s}{\epsilon_{D}^G}\leq FI(\mathbf{X},\mathcal{D})\leq \frac{d}{\epsilon_{D}^G}.$
\item[(ii)] $FI(\mathbf{X},\mathcal{D})=\frac{d}{\epsilon_{D}^G}$ if and only if
 the sub-vectors $\mathbf{X}_{I_j}$, $j=1,...s$, only have
 independent r.v.'s.
\item[(iii)] $FI(\mathbf{X},\mathcal{D})=\frac{s}{\epsilon_{D}^G}$ if and only
if  the sub-vectors $\mathbf{X}_{I_j}$, $j=1,...s$,  only have
totally dependent random r.v.'s.
\end{itemize}
\end{pro}
\dem Just observe that $1\leq \epsilon_{I_j}\leq |I_j|$, with the
lower and upper bounds corresponding to, respectively, complete
dependence and independence of r.v.'s within sub-vectors
$\mathbf{X}_{I_j}$, $j=1,...s$. \fdem\\

Next result presents an extremal coefficient for the amount of
dependence between $\mathbf{Y}_{I_j}$, $j=1,...,s$, $\mathbf{Y}\sim
G$, through the FI of $\mathbf{X}\sim F$ in the domain of attraction
of the MEV $G$.
\begin{pro}\label{p6}
Under the conditions of Proposition \ref{p3}, 
we have
$$
G(\mathbf{x})=(G_{I_1}(\mathbf{x}_{I_1}),...,G_{I_s}(\mathbf{x}_{I_s}))^{1/FI(\mathbf{X},\mathcal{D})}.
$$
\end{pro}
\dem Observe that
$$
\begin{array}{rl}
G(\mathbf{x})=& G_1(x)^{\epsilon_D^G}=(G_1(x)^{\sum_{j=1}^s\epsilon_{I_j}^G})^{1/FI(\mathbf{X},\mathcal{D})}\vspace{0.35cm}\\
=&(G_1(x)^{\epsilon_{I_1}^G}...G_1(x)^{\epsilon_{I_s}^G})^{1/FI(\mathbf{X},\mathcal{D})}\vspace{0.35cm}\\
=&(G_{I_1}(\mathbf{x}_{I_1}),...,G_{I_s}(\mathbf{x}_{I_s}))^{1/FI(\mathbf{X},\mathcal{D})}.\,\,\square
\end{array}
$$

If we consider partition $\mathcal{D}^*$ in the previous result, we obtain the known relation (Smith, \cite{smith} 1990),
$$
G(\mathbf{x})=(G_1(x)^{d})^{1/FI(\mathbf{X},\mathcal{D}^*)}=(G_1(x)^{d})^{\epsilon^G/d}=G_1(x)^{\epsilon^G}.
$$
Observe also that we can write $FI(\mathbf{Y},\mathcal{D})$ instead
of $FI(\mathbf{X},\mathcal{D})$ with $\mathbf{X}$ in the domain of
attraction of the MEV distribution of $\mathbf{Y}$ since
$\mathbf{Y}$ belongs to the same domain of attraction. We finish
this section with a generalization of Proposition \ref{p3} to the
case of equivalent margins and two illustrative examples.
\begin{pro}\label{p7}
If $\mathbf{X}$ has d.f.\hs $F$ with equivalent marginal
distributions in the sense of (\ref{equivmar}), and belongs to the
domain of attraction of a MEV $G$ with unit Fréchet margins then,
for each $k\in\{1,...,s\}$, 
we have
$$
FI(\mathbf{X},\mathcal{D})=\frac{\sum_{j=1}^s
\log C_G(e^{-\gamma_1\mathds{1}_{I_j}(1)},...,e^{-\gamma_d\mathds{1}_{I_j}(d)})}{\log C_G(e^{-\gamma_1},...,e^{-\gamma_d})}\,.
$$\\
\end{pro}

\begin{ex}
We consider a random vector $\mathbf{X}$ of Example 3.2 in Falk and
Tichy (\cite{falk+tichy1}, 2010), i.e., having components
$X_i=\sum_{k=1}^m\lambda_{ik}Y_k$, where $Y_1,...,Y_m$ are
independent r.v.'s with Pareto($\alpha$) distribution, $\alpha>0$,
and $\lambda_{ij}\geq 0$ such that
$\sum_{k=1}^m\lambda_{ik}^{\alpha}=1$, $i=1,...,d$. Taking for $H$
any of the distributions of the margins of $\mathbf{X}$, the
equivalence condition (\ref{equivmar}) holds with $\gamma_i=1$,
$i=1,...,d$. The distribution of $\mathbf{X}$ belongs to the domain
of attraction of
$$
\begin{array}{c}
G(\mathbf{x})=\exp\big(-\sum_{k=1}^m\vee_{i=1}^d\big(\frac{\lambda_{ik}}{x_i}\big)^{\alpha}\big),\,\mathbf{x}>\mathbf{0}
\end{array}
$$
with Fréchet margins, $G_i(x)=\exp(-x^{-\alpha})$. Hence, for
$\mathbf{u}=(u_1,...,u_d)\in(\mathbf{0},\mathbf{1})$,
$$
\begin{array}{c}
C_G(\mathbf{u})=\exp\big(-\sum_{k=1}^m\vee_{i=1}^d\lambda_{ik}^{\alpha}(-\log u_i)\big).
\end{array}
$$
By Proposition \ref{p7}, we have
$$
\begin{array}{c}
\dst FI(\mathbf{X},\mathcal{D}^*)=\frac{\sum_{j=1}^d  \sum_{k=1}^m\vee_{i=1}^d \lambda_{ik}^{\alpha}\mathds{1}_{\{j\}}(i)
}{ \sum_{k=1}^m\vee_{i=1}^d \lambda_{ik}^{\alpha}}=\frac{\sum_{j=1}^d  \sum_{k=1}^m \lambda_{jk}^{\alpha}
}{ \sum_{k=1}^m\vee_{i=1}^d \lambda_{ik}^{\alpha}}=\frac{d}{\sum_{k=1}^m\vee_{i=1}^d \lambda_{ik}^{\alpha}}\,.
\end{array}
$$
as obtained in Falk and Tichy (\cite{falk+tichy1}, 2010). For any
partition $\mathcal{D}$, we have
$$
\begin{array}{c}
\dst FI(\mathbf{X},\mathcal{D})=\dst\frac{\sum_{j=1}^s  \sum_{k=1}^m\vee_{i=1}^d \lambda_{ik}^{\alpha}\mathds{1}_{I_j}(i)
}{ \sum_{k=1}^m\vee_{i=1}^d \lambda_{ik}^{\alpha}}
\,.
\end{array}
$$
To illustrate, consider $d=3=m$, $\alpha =1$ and weights
$$
\begin{array}{ccc}
\lambda_{11}=4/8,&  \lambda_{12}=2/8,&\lambda_{13}=2/8,\\
\lambda_{21}=1/8,&  \lambda_{22}=1/8,&\lambda_{23}=6/8,\\
\lambda_{31}=3/8,&  \lambda_{32}=2/8,&\lambda_{33}=3/8.
\end{array}
$$
We have
$$
\begin{array}{rl}
FI(\mathbf{X},\mathcal{D}^*)=&\frac{3}{4/8+2/8+6/8}=\frac{24}{12}=2
\end{array}
$$
and, for $\mathcal{D}=\{\{1,2\},\{3\}\}$,
$$
\begin{array}{rl}
FI(\mathbf{X},\mathcal{D})=&\frac{(4/8+2/8+6/8)+(3/8+2/8+3/8)}{4/8+2/8+6/8}=\frac{20}{12}<2.
\end{array}
$$
\end{ex}

\begin{ex}
If $G$ has copula
$$
\begin{array}{c}
C_G(u_1,...,u_d)=\exp\big(-\big(\sum_{i=1}^d(-\log u_i)^{1/\alpha}\big)^{\alpha}\big),\,0<\alpha\leq 1,
\end{array}
$$
(symmetric logistic model) then, for any partition $\mathcal{D}$, we
have
$$
\begin{array}{c}
FI(\mathbf{X},\mathcal{D})=\frac{\sum_{j=1}^s|I_j|^{\alpha}}{d^{\alpha}}=\sum_{j=1}^s\big(\frac{|I_j|}{d}\big)^{\alpha},
\end{array}
$$
and $FI(\mathbf{X},\mathcal{D}^*)=d^{1-\alpha}$ as already stated in
Geluk \emph{et al.} (\cite{gel+}, 2007). In the symmetric model the
FI is only a function of the blocks size. If we consider the more
general asymmetric logistic model, whose copula is given by
$$
\begin{array}{c}
C_G(u_1,...,u_d)=\exp\big\{-\sum_{k=1}^q\big(\sum_{i=1}^d(-\beta_{ki}\log u_i)^{1/\alpha_k}\big)^{\alpha_k}\big\}
\end{array}
$$
where $\beta_{ki}$ are non-negative constantes such that
$\sum_{k=1}^q\beta_{ki}=1$, $i=1,...,d$, $0<\alpha_k\leq 1$,
$k=1,...,q$, we obtain
$$
\begin{array}{c}
FI(\mathbf{X},\mathcal{D})=\frac{\sum_{j=1}^s\sum_{k=1}^q\big(\sum_{i\in I_j}\beta_{ki}^{1/\alpha_k}\big)^{\alpha_k}}
{\sum_{k=1}^q\big(\sum_{i=1}^d\beta_{ki}^{1/\alpha_k}\big)^{\alpha_k}}.
\end{array}
$$
\end{ex}

\section{Tail dependence for blocks}\label{sli}

{In the following we always consider that $\mathbf{X}$ has
continuous marginal d.f.'s.} Consider notation $M(I_j)=\bigvee_{i\in
I_j}F_i(X_i)$, $j\in D$.
\begin{defi}
The upper-tail dependence coefficients of { $\mathbf{X}$}
corresponding to partition $\mathcal{D}$ of $D$ are defined by, for
each $S\subsetneq\{1,...,s\}$,
\begin{eqnarray}\label{tau}
\begin{array}{c}
{\tau_S^F}=\lim_{u\uparrow 1}P(\bigcap_{j\not\in S}M(I_j)>u|\bigcap_{j\in S}M(I_j)>u)
\end{array}
\end{eqnarray}
when the limit exists.
\end{defi}
If we consider partition $\mathcal{D}^*$, then $S\subsetneq\{1,...,d\}$ and we find the definition of Li (\cite{li3}, 2009).
Further, the case $s=2$ lead us to definition of Ferreira and Ferreira (\cite{hf+mf3}, 2011).\\

{Consider 
\begin{eqnarray}\label{lambda}
{\lambda_S^F}:=\lim_{u\uparrow 1}\frac{P(\bigcap_{j\in S}M(I_j)>u)}{1-u}
\end{eqnarray}
for each $S\subsetneq\{1,...,s\}$}. Hence we can write
{\begin{eqnarray}\label{taulambda}
\tau_S^F=\frac{\lambda_{\{1,...,s\}}^F}{\lambda_S^F}.
\end{eqnarray}
Observe that $\lambda_S^F$ corresponds to the multivariate
upper-tail dependence coefficient $\Lambda_U(\mathbf{1}_S)$
 in Schmidt and Stadtm\"{u}ller
(\cite{schmidt+stad}, 2006), where $\mathbf{1}_S$ denotes the unit
vector with dimension $|S|$. In particular, for partition
$\mathcal{D}^*$, $\lambda_{\{i,j\}}^F=\Lambda_U(1,1)$
 corresponds to the well-known bivariate tail dependence concept (Sibuya, \cite{sib} 1960; Joe, \cite{joe} 1993).\\
}

Before we relate the FI with the tail dependence coefficients corresponding to a partition, we present in this section
some extensions of the results in Li (\cite{li3}, 2009).

\begin{pro}\label{p8}
If $\mathbf{X}$ has MEV distribution $G$ with standard Fréchet
margins and spectral measure $W$ defined on the $d$-dimensional unit
sphere $S_d$ then, for each $S\subsetneq\{1,...,s\}$, we have
\begin{eqnarray}\label{taulambda}
\tau_S^G=\frac{\int_{S_d}\bigwedge_{j=1}^s\bigvee_{i\in I_j}w_i\,dW(\mathbf{w})}{\int_{S_d}\bigwedge_{j\in s}\bigvee_{i\in I_j}w_i\,dW(\mathbf{w})}.
\end{eqnarray}
\end{pro}
\dem From the spectral representation of $G$ we obtain, for $u$
sufficiently close to $1$,
\begin{eqnarray}\label{p8.1}
\begin{array}{rl}
P\big(\bigcap_{j=1}^s M(I_j)>u\big)=&1-\sum_{\emptyset\not = S\subset \{1,...,s\}}(-1)^{|S|+1}
G\Big(-\frac{\mathds{1}_{\mathcal{I}(S)}^{-1}(1)}{\log u},...,-\frac{\mathds{1}_{\mathcal{I}(S)}^{-1}(d)}{\log u}\Big)\vspace{0.35cm}\\
\approx& 1-\sum_{\emptyset\not = S\subset \{1,...,s\}}(-1)^{|S|+1}\big(1+\log u\int_{S_d}\bigvee_{i\in \mathcal{I}(S)}w_i\,dW(\mathbf{w})\big)\vspace{0.35cm}\\
=& 1-\sum_{\emptyset\not = S\subset \{1,...,s\}}(-1)^{|S|+1}\big(1-(-\log u)\int_{S_d}\bigvee_{j\in S}\bigvee_{i\in I_j}w_i\,dW(\mathbf{w})\big).
\end{array}
\end{eqnarray}
Since
\begin{eqnarray}\nn
\begin{array}{c}
  \sum_{\emptyset\not = S\subset \{1,...,s\}}(-1)^{|S|+1}=1
\end{array}
\end{eqnarray}
and
\begin{eqnarray}\nn
\begin{array}{c}
  \sum_{\emptyset\not = S\subset \{1,...,s\}}(-1)^{|S|+1}\bigvee_{j\in S}a_j=\bigwedge_{j\in\{1,...,s\}}a_j,
\end{array}
\end{eqnarray}
expression in (\ref{p8.1}) becomes
\begin{eqnarray}\nn
\begin{array}{c}
P\big(\bigcap_{j=1}^s M(I_j)>u\big)\approx(-\log u)\int_{S_d}\bigwedge_{j=1}^s\bigvee_{i\in I_j}w_i\,dW(\mathbf{w})\big).
\end{array}
\end{eqnarray}
Analogously, we obtain
\begin{eqnarray}\nn
\begin{array}{c}
P\big(\bigcap_{j\in S}M(I_j)>u\big)\approx(-\log u)\int_{S_d}\bigwedge_{j\in S}\bigvee_{i\in I_j}w_i\,dW(\mathbf{w})\big).\,\,\square
\end{array}
\end{eqnarray}
For the particular case $\mathcal{D}^*$, the previous result is the
one found in Li (\cite{li3}, 2009). Note also that the numerator of
(\ref{taulambda}) can be expressed through extremal coefficients as
follows:
\begin{eqnarray}\nn
\begin{array}{rl}
&\sum_{\emptyset\not = S\subset \{1,...,s\}}(-1)^{|S|+1}\int_{S_d}\bigvee_{j\in S}\bigvee_{i\in I_j}w_i\,dW(\mathbf{w})\vspace{0.35cm}\\
=&\sum_{\emptyset\not = S\subset \{1,...,s\}}(-1)^{|S|+1}\epsilon_{\mathcal{I}(S)}^G\vspace{0.35cm}\\
=&\epsilon_{I_1}^G+...+\epsilon_{I_s}^G-(\epsilon_{I_1\cup I_2}^G+\epsilon_{I_1\cup I_3}^G+...)-...+(-1)^{|S|+1}\epsilon_{I_1\cup...\cup I_s}^G,
\end{array}
\end{eqnarray}
a generalization of result (15) in Ferreira and Ferreira
(\cite{hf+mf3}, 2011) where $s=2$.\\

The next result highlights the connections between tail dependence
and extremal coefficients.
\begin{cor}
Under the conditions of Proposition \ref{p8}, we have
\begin{itemize}
\item[(i)] $\lambda_S^G=\sum_{\emptyset\not = T\subset
S}(-1)^{|T|+1}\epsilon_{\mathcal{I}(T)}^G$
\item[(ii)]$\tau_S^G=\frac{\sum_{\emptyset\not = T\subset
\{1,...,s\}}(-1)^{|T|+1}\epsilon_{\mathcal{I}(T)}^G}{\sum_{\emptyset\not
= T\subset S}(-1)^{|T|+1}\epsilon_{\mathcal{I}(T)}^G}$.
\end{itemize}
\end{cor}

We end this section with a generalization of Theorem 2.6 in Li
(\cite{li3}, 2009), by adapting the arguments to subsets of $D$ that
correspond to unions of blocks in $\mathcal{D}$.

\begin{pro}\label{p9}
If $F${ belongs to the domain of attraction of a MEV $G$ with unit
Fréchet margins} then, for any partition $\mathcal{D}$ and
$\emptyset\not=S\subset \{1,...,s\}$, the non-null upper-tail
dependence coefficients $\tau_S^F$ are the same as the corresponding
ones of $G$.
\end{pro}

\section{Fragility Index and tail dependence for
blocks}\label{sli+fi}
 In this section we shall see that the asymptotic d.f.'s of the
 conditional probability of $k$ exceedances between blocks,
 $I_1,...,I_s$, can be derived
through the tail dependence coefficients given in (\ref{lambda}).
More precisely, if the d.f.\hs $F$ of $\mathbf{X}$ has tail
dependence coefficients $\lambda_S^F$ corresponding to partition
$\mathcal{D}$, we can obtain $\lim_{x\to\infty}P(N_x=k|N_x>0)$ from
these latter. In particular, for $F$ in the domain of attraction of
a MEV $G$, besides the representations presented in the previous
results, we can write those limiting probabilities through tail
dependence coefficients $\lambda_S^G$.

\begin{pro}\label{p10}
Let $\mathbf{X}$ be a random vector with d.f.\hs $F$ with continuous
and identically distributed margins. Let $\mathcal{D}$ be a
partition of $D$ for which the tail dependence coefficients
$\lambda_S^F$ corresponding to $F$ exist for each $S\subset
\{1,...,s\}$. Then,

\begin{itemize}
\item[(i)] for each $k\in\{1,...,s\}$,

$$\lim_{x\to\infty}P(N_x=k|N_x>0)=\frac{\sum_{S\subset
\{1,...,s\},|S|=k}\sum_{T\subset S^C}(-1)^{|T|}\lambda_{T\cup
S}^F} {\sum_{k=1}^s\sum_{S\subset
\{1,...,s\},|S|=k}\sum_{T\subset S^C}(-1)^{|T|}\lambda_{T\cup
S}^F},$$

as long as the numerator is non-null.
\item[(ii)]if $F$ belongs to the domain of attraction of a MEV $G$ with unit
Fréchet marginals, the limits in (i) exist and coincide for both
distributions, $F$ and $G$.
\end{itemize}
\end{pro}
\dem Just observe that
$$
\begin{array}{c}
 \lim_{x\to\infty}P(N_x=k|N_x>0)=\dst\lim_{u\uparrow 1}\frac{\sum_{S\subset
\{1,...,s\},|S|=k}\sum_{T\subset S^C}(-1)^{|T|}P(\bigcap_{T\cup
S}M(I_j)>u)} {\sum_{k=1}^s\sum_{S\subset
\{1,...,s\},|S|=k}\sum_{T\subset S^C}(-1)^{|T|}P(\bigcap_{T\cup
S}M(I_j)>u)}.
\end{array}
$$
Now divide both terms of the ratio by $1-u$. \fdem\\

Analogously, we obtain the FI through $\lambda_S^F$ or
$\lambda_S^G$.
\begin{pro}\label{p11}
Under the conditions of Proposition \ref{p10}, we have
\begin{itemize}
\item[(i)]
$$FI(\mathbf{X},\mathcal{D})=\frac{\sum_{j=1}
^s\lambda_{\{j\}}^F} {\sum_{k=1}^s\sum_{S\subset
\{1,...,s\},|S|=k}\sum_{T\subset S^C}(-1)^{|T|}\lambda_{T\cup
S}^F}$$
as long as the numerator is non-null.
\item[(ii)]if $F$ belongs to the domain of attraction of a MEV $G$ with unit
Fréchet marginals, $FI(\mathbf{X},\mathcal{D})$ exists and the
the expression in (i) coincides for both distributions, $F$ and
$G$.
\end{itemize}
\end{pro}
The statements in (ii) of the two propositions are consequences of
Proposition \ref{p9}.

\section{Asymptotic Independence}\label{saifi}
If $\mathbf{X}$ has independent margins $X_{i}$, $i=1,...,d$, we
have an unit FI. As Geluk \emph{et al.} (\cite{gel+}, 2007)
observed, in this case we might have an asymptotic independence
characterized by a dependency that vanishes at extreme quantiles.
Asymptotic independence means that the system is weakly fragile,
albeit possibly correlated (e.g., gaussian vectors). Geluk \emph{et
al.} (\cite{gel+}, 2007) have defined a fragility index for
asymptotic independence (AIFI) given by
\begin{eqnarray}\label{aifiG}
\eta_D=\frac{1}{d}\lim_{x\to\infty}\frac{\sum_{i=1}^{d}
\log P(X_i>x)}{\log P(X_1>x,...,X_d>x)}.
\end{eqnarray}
In case $d=2$  we find  the Ledford and Tawn coefficient of
asymptotic independence (\cite{led+tawn1,led+tawn2}, 1996, 1997)
and,  if $d>2$, a multivariate extension of this latter (Ferreira
and
Ferreira, \cite{hf+mf2}, 2012).\\


%
%
%

 Here we consider an extension of the AIFI in (\ref{aifiG}) for
 blocks, in the same spirit of the FI in Proposition \ref{p3},
 i.e., by relating the AIFI within blocks with the AIFI of the whole vector.

Let $\eta_A$ be the AIFI of sub-vector $\mathbf{X}_A$ of
$\mathbf{X}$, with $A\subset D$, i.e.,
\begin{eqnarray}\label{aifiA}
\eta_A=\frac{1}{|A|}\lim_{x\to\infty}\frac{\sum_{i\in A}
\log P(X_i>x)}{\log P(\mathbf{X}_A>\mathbf{x}_A)}.
\end{eqnarray}

\begin{defi}
Let $\mathbf{X}=(X_1,...,X_d)$ be a random vector with
$FI(\mathbf{X},\mathcal{D})=1$. Then the AIFI of
$\mathbf{X}=(X_1,...,X_d)$ relative to partition $\mathcal{D}$ of
$D$ is
\begin{eqnarray}\label{aifi}
\frac{1}{s}\lim_{x\to\infty}\frac{\sum_{j=1}^s
\log P(\mathbf{X}_{I_j}>\mathbf{x}_{I_j})}{\log P(X_1>x,...,X_d>x)},
\end{eqnarray}
whenever the limit exists and is denoted
$\eta(\mathbf{X},\mathcal{D})$.
\end{defi}
\begin{pro}\label{pind1}
Let $\mathbf{X}=(X_1,...,X_d)$ be a random vector with
$FI(\mathbf{X},\mathcal{D})=1$. Assume that (\ref{aifiG}) holds and
that (\ref{aifiA}) holds for all $I_j\in \mathcal{D}$, $j=1,...,s$,
with limit given by $\eta_{I_j}$, respectively. If $\mathbf{X}$ has
identically distributed or equivalent margins in the sense of
(\ref{equivmar}), then
\begin{eqnarray}\label{pind1.1}
\eta(\mathbf{X},\mathcal{D})=\eta_D\frac{1}{s}\sum_{j=1}^{s}\frac{1}{\eta_{I_j}}.
\end{eqnarray}

\dem Observe that by (\ref{aifiG}) and (\ref{aifiA}), then
\begin{eqnarray}\nn
\eta(\mathbf{X},\mathcal{D})=\lim_{x\to\infty}\eta_D\,\,
\frac{\frac{1}{s}\sum_{j=1}^{s}\frac{1}{\eta_{I_j}|I_j|}\sum_{i\in I_j}\log P(X_i>x)}{\frac{1}{d}\sum_{i=1}^{d}\log P(X_i>x)},
\end{eqnarray}
and by (\ref{equivmar}), we have
\begin{eqnarray}\nn
\eta(\mathbf{X},\mathcal{D})=\lim_{x\to\infty}\eta_D\,\,
\frac{\frac{1}{s}\sum_{j=1}^{s}\frac{1}{\eta_{I_j}}\big(\frac{1}{|I_j|}\sum_{i\in I_j}\log \gamma_i+\log(1-H(x))\big)}
{\frac{1}{d}\sum_{i=1}^{d}\log \gamma_i+\log(1-H(x))}\,.\,\,\square  
\end{eqnarray}
 \end{pro}

In particular, we have $\eta(\mathbf{X},\mathcal{D}^*)=\eta_D$ and
hence, the AIFI of the system divided into blocks is larger than the
AIFI of the system itself, i.e., $\eta(\mathbf{X},\mathcal{D})\geq
\eta(\mathbf{X},\mathcal{D}^*)$.

\begin{rem}
In the particular case of $FI(\mathbf{X},\mathcal{D}^*)>1$ we have
$\eta(\mathbf{X},\mathcal{D})=1$.

\nid Observe also that, by Proposition \ref{pind1}, we can relate
$\eta(\mathbf{X},\mathcal{D})$ with the fragility indexes of the
whole system and of each block. More precisely,
$\eta(\mathbf{X},\mathcal{D})$ is the arithmetic mean of the ratios
$\eta(\mathbf{X},\mathcal{D}^*)/\eta(\mathbf{X}_{I_j},\mathcal{D}^*)$.


\end{rem}
%

\begin{ex}
Consider $\mathbf{X}=(X_1,...,X_d)$ a standard $d$-variate Gaussian
random vector with d.f.\hs $\Phi_d( \cdot;(\Sigma_{i,j})_{i,j\in
D})$ having positive definite correlation matrix
$(\Sigma_{i,j})_{i,j\in D}$. We have, for all $A\subset D$,
$\eta_A^{-1}= \mathbf{1}_A(\Sigma_{i,j})^{-1}_{i,j\in
A}\mathbf{1}_A^T$, the sum of all elements of the sub-matrix
$(\Sigma_{i,j})^{-1}_{i,j\in A}$ (Geluk \emph{et al.}, \cite{gel+},
2007; Hua and Joe, \cite{hua+joe}, 2011). For illustration, consider
dimension $d=4$, constant correlation $\rho$, and take $s=3$ with
$I_1=\{1,2\}$, $I_2=\{3\}$ and $I_3=\{4\}$. We have $\eta_D=(3\rho
+1)/4$, $\eta_{I_1}=(\rho +1)/2$ and $\eta_{I_2}=\eta_{I_3}=1$.
Hence
$$\begin{array}{l}
\eta(\mathbf{X},\mathcal{D})=
\frac{1}{3}\frac{3\rho +1}{4}\big(\frac{2}{\rho +1}+1+1\big)
=\frac{(3\rho +1)(\rho +2)}{6(\rho +1)}.\\
\end{array}
$$
\end{ex}


In the sequel we consider positive/negative association of a random
vector in the sense of Ledford and Tawn (\cite{led+tawn1,led+tawn2},
1996, 1997).

\begin{pro}[inter-blocks asymptotic independence]\label{pind2}
Under the conditions of Proposition \ref{pind1}, 
we have
\begin{itemize}
\item[(i)] $\eta(\mathbf{X},\mathcal{D})\leq
\frac{1}{s}$ in case of positive association between sub-vectors
$\mathbf{X}_{I_j}$, $j=1,...s$, and
$\eta(\mathbf{X},\mathcal{D})\geq \frac{1}{s}$ in case of
negative association.
\item[(ii)] $\eta(\mathbf{X},\mathcal{D})\geq\frac{\wedge_{j=1}^s\eta_{I_j}}{s}\sum_{j=1}^{s}\frac{1}{\eta_{I_j}}$.
\item[(iii)] $\eta(\mathbf{X},\mathcal{D})=\frac{1}{s}$ if and only if
 the sub-vectors $\mathbf{X}_{I_j}$, $j=1,...s$, are
 independent.
 \item[(iv)] $\eta(\mathbf{X},\mathcal{D})=\frac{\wedge_{j=1}^s\eta_{I_j}}{s}\sum_{j=1}^{s}\frac{1}{\eta_{I_j}}$ if and only
if  the sub-vectors $\mathbf{X}_{I_j}$, $j=1,...s$,  are totally
dependent.
\end{itemize}
\end{pro}
\dem Observe that, $P(X_1>x,...,X_d>x)\leq
\wedge_{j=1}^sP(\mathbf{X}_{I_j}>\mathbf{x}_{I_j}) $, with the upper
bound corresponding to total dependence between sub-vectors
$\mathbf{X}_{I_j}$, $j=1,...s$. Under positive association we have
$P(X_1>x,...,X_d>x)\geq
\prod_{j=1}^sP(\mathbf{X}_{I_j}>\mathbf{x}_{I_j})$ and negative
association otherwise, with the bounds corresponding to independence
between sub-vectors
$\mathbf{X}_{I_j}$, $j=1,...s$. \fdem\\

\begin{pro}[intra-blocks asymptotic independence]\label{pind3}
Under the conditions of Proposition \ref{pind1}, 
we have
\begin{itemize}
\item[(i)] $\eta(\mathbf{X},\mathcal{D})\leq
\frac{\eta_Dd}{s}$ in case the sub-vectors $\mathbf{X}_{I_j}$,
$j=1,...s$, only have positively associated r.v.'s and
$\eta(\mathbf{X},\mathcal{D})\geq \frac{\eta_Dd}{s}$ in case of
negative association.
\item[(ii)] $\eta(\mathbf{X},\mathcal{D})\geq\eta_{D}$.
\item[(iii)] $\eta(\mathbf{X},\mathcal{D})=\frac{\eta_Dd}{s}$ if and only if
 the sub-vectors $\mathbf{X}_{I_j}$, $j=1,...s$, only have
 independent r.v.'s.
 \item[(iv)] $\eta(\mathbf{X},\mathcal{D})=\eta_{D}$ if and only
if  the sub-vectors $\mathbf{X}_{I_j}$, $j=1,...s$,  only have
totally dependent random r.v.'s.
\end{itemize}
\end{pro}
\dem For each  $j=1,...s$, under positive association of the  r.v.'s
in $\mathbf{X}_{I_j}$ we have
$P(\mathbf{X}_{I_j}>\mathbf{x}_{I_j})\geq \prod _{i\in
I_j}P(X_i>x)$, and hence $ 1/ \eta_{I_j}\leq |I_j|$, with the upper
bound corresponding to
independence. For negative association we have $P(\mathbf{X}_{I_j}>\mathbf{x}_{I_j})\leq \prod _{i\in I_j}P(X_i>x)$.
Observe also that total dependence within each block means $P(\mathbf{X}_{I_j}>\mathbf{x}_{I_j})=\wedge_{i\in I_j}\gamma_i(1-H(x))$. \fdem\\

As already mentioned, the definition of the AIFI in (\ref{aifi})
measures the asymptotic independent fragility of a system divided
into blocks (sub-vectors), by relating the asymptotic independent
fragility within the blocks and in the whole system. If in
(\ref{aifiG}) we generalize for blocks the concept of an exceedance
of a r.v., $X_i>x$, through events $\mathbf{X}_{I_j}\not\leq
\mathbf{x}_{I_j}$, we obtain another coefficient for asymptotic tail
independence. In this way, we extend the \emph{coefficient of
asymptotic tail independence} that was considered in Ferreira and
Ferreira (\cite{hf+mf3}, 2011) for the particular case of a
partition $\mathcal{D}=\{I_1,I_2\}$ of $D=\{1,...,d\}$.

\begin{defi}
Let $\mathbf{X}=(X_1,...,X_d)$ be a random vector with
$FI(\mathbf{X},\mathcal{D})=1$. The coefficient of asymptotic
independence of $\mathbf{X}=(X_1,...,X_d)$ relative to partition
$\mathcal{D}$ of $D$
is
\begin{eqnarray}\label{aifi2}
\lim_{x\to\infty}\frac{1}{s}\frac{\sum_{j=1}^s\log P(\mathbf{X}_{I_j}\not\leq\mathbf{x}_{I_j})}
{\log P(\bigcap_{j=1}^s\{\mathbf{X}_{I_j}\not\leq\mathbf{x}_{I_j}\})},
\end{eqnarray}
whenever the limit exists, and is denoted $\eta_{(I_1,...,I_s)}$.
\end{defi}
The following result is therefore an immediate extension of
Proposition 2.4 in Ferreira and Ferreira (\cite{hf+mf3}, 2011).

\begin{pro}
Let $\mathbf{X}=(X_1,...,X_d)$ be a random vector with
$FI(\mathbf{X},\mathcal{D})=1$. Assume that the limit in
(\ref{aifi2}) exists and,
 for all $\emptyset\not=K_j\subset I_j$, $j=1,...,s$, (\ref{aifiA}) holds for $\cup_{j=1}^sK_j$. Then
\begin{eqnarray}\label{p2aifi2}
\eta_{(I_1,...,I_s)}=\max\{\eta_{\{i_1,...,i_s\}}:\,i_j\in
I_j,\,j=1,...,s\}.
\end{eqnarray}
\end{pro}

\nid Similar to $\eta(\mathbf{X},\mathcal{D})$, coefficient
$\eta_{(I_1,...,I_s)}$ is also based on the coefficient
of Ledford and Tawn or multivariate extensions of this latter.\\

In the example below, one can see that the asymptotic tail
independent coefficients, $\eta(\mathbf{X},\mathcal{D})$ and
$\eta_{(I_1,...,I_s)}$, are different.

\begin{ex}
Consider $\{V_i\}_{i\geq 1}$ an i.i.d.\hs sequence of unit Pareto
r.v.'s. Let $(X_1,X_2,X_3)$ be a random vector such that,
$X_1=\min(V_1,V_2)$, $X_2=\min(V_2,V_3)$ and $X_3=\min(V_3,V_4)$. If
$\mathcal{D}=\{I_1,I_2\}$, with $I_1=\{1,2\}$ and $I_2=\{3\}$, we
have $P(X_i>x)=x^{-2}$, $i=1,2,3$,
$P(X_1>x,X_2>x)=P(X_2>x,X_3>x)=x^{-3}$ and
$P(X_1>x,X_3>x)=P(X_1>x,X_2>x,X_3>x)=x^{-4}$. Therefore
$$
FI(\mathbf{X},\mathcal{D})=\lim_{x\to\infty}=\frac{P(\bigcup_{i=1}^2 \{X_i>x\})+P(X_3>x)}{ P(\bigcup_{i=1}^3\{X_i>x\})}=\frac{2x^{-2}-x^{-3}+x^{-2}}{3x^{-2}-2x^{-3}-x^{-4}+x^{-4}}=1,
$$
$$\begin{array}{l}\eta(\mathbf{X},\mathcal{D})=\frac{1}{2}\lim_{x\to\infty}\frac{\sum_{j=1}^2\log
P(\mathbf{X}_{I_j}>\mathbf{x}_{I_j})}{\log
P(X_1>x,X_2>x,X_3>x)}=\frac{1}{2}\lim_{x\to\infty}\frac{\log
P(X_1>x,X_2>x)+P(X_3>x)}{\log
P(X_1>x,X_2>x,X_3>x)}=\frac{1}{2}\frac{-3-2}{-4}=\frac{5}{8}\end{array}$$
and
$$\begin{array}{l}\eta_{(I_1,I_2)}=\frac{1}{2}\lim_{x\to\infty}\frac{\log P(\mathbf{X}_{I_1}\not\leq \mathbf{x}_{I_1})+\log P(\mathbf{X}_{I_2}\not\leq \mathbf{x}_{I_2})}{\log P(\mathbf{X}_{I_1}\not\leq \mathbf{x}_{I_1},\mathbf{X}_{I_2}\not\leq \mathbf{x}_{I_2})}\vspace{0.35cm}\\
=\frac{1}{2}\lim_{x\to\infty}\frac{\log (P(X_1>x)+P(X_2>x)-P(X_1>x,X_2>x))+\log P(X_3>x)}{\log (P(X_1>x,X_3>x)+P(X_2>x,X_3>x)-P(X_1>x,X_2>x,X_3>x))}\vspace{0.35cm}\\
=\frac{1}{2}\lim_{x\to\infty}\frac{\log (2x^{-2}-x^{-3})+\log x^{-2}}{\log (x^{-4}+x^{-3}-x^{-4})}=2/3.\end{array}$$\\

Observe that the same results are obtained if we apply Propositions
6.1 and 6.4, respectively. More precisely, we have
$$
\eta_D=\eta_{\{1,2,3\}}=\frac{1}{3}\lim_{x\to\infty}\frac{\sum_{i=1}^3\log P(X_i>x)}{\log P(X_1>x,X_2>x,X_3>x)}=\frac{1}{3}\frac{3\log x^{-2}}{\log x^{-4}}=\frac{1}{2},
$$
$$
\begin{array}{ccc}\dst\eta_{\{1,3\}}=\frac{1}{2}\frac{2\log x^{-2}}{\log x^{-4}}=\frac{1}{2}
&\textrm{and}&
\dst\eta_{\{2,3\}}=\frac{1}{2}\frac{2\log x^{-2}}{\log x^{-3}}=\frac{2}{3}.
\end{array}$$
Obviously $\eta_{I_2}=1$ and
$$
\eta_{I_1}=\frac{1}{2}\lim_{x\to\infty}\frac{\sum_{i=1}^2\log P(X_i>x)}{\log P(X_1>x,X_2>x)}=\frac{1}{2}\frac{2\log x^{-2}}{\log x^{-3}}=\frac{2}{3}.
$$

Hence, by Propositions 6.1 and 6.4, respectively,
$$\eta(\mathbf{X},\mathcal{D})=\frac{1}{2}\,\Big(\frac{1}{2}(3/2+1)\Big)=\frac{5}{8}$$
and
$$\eta_{(I_1,I_2)}=max\{\eta_{\{1,3\}},\eta_{\{2,3\}}\}=max\{1/2,2/3\}=2/3.$$
\end{ex}

\section{Estimation of the Fragility Index for blocks}\label{sestim}
In the previous sections we have relate the FI for blocks with well-known tail dependence measures. This will allow to obtain immediate estimators for our index through the estimators of those coefficients that are already studied in literature.


Proposition \ref{p3} presents an estimation procedure for the FI
based on the extremal coefficients of Tiago de Oliveira
(\cite{tiago}, 1962/63) and Smith (\cite{smith}, 1990) given in
(\ref{extremalcoef}). Observe that they can be expressed  through
the stable tail dependence function, $l_G$, defined in (\ref{stdf}).
 There are
several references in literature on the estimation of the stable
tail dependence function. In a parametric framework, a model for
$l_G$ must be imposed. Non-parametric estimators
 are usually based on an arbitrarily chosen
parameter, corresponding to the number of top order statistics to
use in order to provide the best trade-off between bias and
variance, which is not an easy task. For a survey, we refer Krajina
(\cite{krajina}, 2010) or Beirlant \emph{et al.\hs} (\cite{beirl+}
2004). The more recent work in Ferreira and Ferreira (\cite{hf+mf3},
2011) presents a simpler non-parametric estimator based on sample
means, which we shall adopt here.

For $A\subset D$, denote $M(A)=\bigvee_{i\in A}F_i(X_i)$. Consider
\begin{eqnarray}\label{estimcoefext}
\begin{array}{ccc}
\widehat{\epsilon}_D^G=\frac{\overline{M(D)}}{1-\overline{M(D)}}&\textrm{ and } &
\widehat{\epsilon}_{I_j}^G=\frac{\overline{M(I_j)}}{1-\overline{M(I_j)}},
\end{array}
\end{eqnarray}
where $\overline{M(A)}$ is the sample mean
\begin{eqnarray}\label{estimean}
\begin{array}{c}
\dst\overline{M(A)}=\frac{1}{n}\sum_{i=1}^n\bigvee_{j\in A}\widehat{F}_j(X^{(i)}_j)
\end{array}
\end{eqnarray}
and $\widehat{F}_j$, $j\in A$, is the (modified) empirical d.f. of
$F_j$,
$$
\widehat{F}_j(u)=\frac{1}{n+1}\sum_{k=1}^n\mathds{1}_{\{X_j^{(k)}\leq
u\}}.
$$
The denominator $n+1$ instead of $n$ in the ordinary empirical d.f.
concerns estimation accuracy and other modifications can be used
(see, for instance, Beirlant \emph{et al.\hs} \cite{beirl+} 2004).
Based on Proposition \ref{p3}, and for a partition $\mathcal{D}$ of
$D$, we consider estimator
\begin{eqnarray}\label{estim1}
\begin{array}{ccc}
\widehat{FI}(\mathbf{X},\mathcal{D})=\frac{\sum_{j=1}^s\widehat{\epsilon}_{I_j}^G}{\widehat{\epsilon}_D^G},
\end{array}
\end{eqnarray}
which is consistent given the consistency of estimators
$\widehat{\epsilon}_D^G$ and $\widehat{\epsilon}_{I_j}^G$ already
stated in Ferreira and Ferreira (\cite{hf+mf3}, 2011).

By Proposition \ref{p11}, we can also estimate
$FI(\mathbf{X},\mathcal{D})$ based on the tail dependence
coefficients in (\ref{lambda}). As already mentioned, these
correspond to multivariate upper-tail dependence coefficients
considered in Schmidt and Stadtm\"{u}ller (\cite{schmidt+stad},
2006), for which non-parametric estimators have been studied.
We remark that these estimators are also based on a similar
procedure as described above for the stable tail dependence
function, i.e., it comprises the choice of a number of top order
statistics to be used, not too large neither too small to avoid,
respectively, large bias and large variance.\\


\begin{rem}
In case of asymptotic independence considered in Section
\ref{saifi}, immediate estimators for the asymptotic independent
coefficients, $\eta(\mathbf{X},\mathcal{D})$ and
$\eta_{(I_1,...,I_s)}$, can be derived from the found relations
 with the
Ledford and Tawn coefficient (or multivariate extensions), whose
estimation has already been  studied in literature (see, for
instance, Draisma \cite{draisma+}, 2004 or, for a survey, Beirlant
\emph{et al.}, \cite{beirl+}, 2004). In particular, the Ledford and
Tawn coefficient can be estimated as the extreme value index of
r.v.\hs $\min(1/(1-F_1(X_1)),1/(1-F_2(X_2)))$ and, similarly, its
multivariate versions in (\ref{aifiA}) can be estimated as the
extreme value index of r.v.\hs $\min_{i\in A}(1/(1-F_i(X_i)))$.\\
\end{rem}

\subsection{An application to financial data}\label{saplic}

We now illustrate the estimation of the FI for blocks through an
application to data analyzed by Ferreira and Ferreira
(\cite{hf+mf3}, 2011). The data are the series of negative
log-returns of the closing values of the stock market indexes, CAC
40 (France), FTSE100 (UK), SMI (Swiss), XDAX (German), Dow Jones
(USA), Nasdaq (USA), SP500 (USA), HSI (China), Nikkei (Japan). The
period covered is January 1993 to March 2004. Since we do not have a
sample of maximum values, we consider the monthly maximums in each
market. We group the indexes in Europe (CAC 40, FTSE100, SMI, XDAX),
USA (Dow Jones, Nasdaq) and Far East (HSI, Nikkei). The presence of
dependence within these groups was already evidenced in Ferreira and
Ferreira (\cite{hf+mf3}, 2011). We are interested in assessing the
fragility within the system of the financial stock markets whenever
grouped in the three big world markets referred: Europe, USA and Far
East. To this end, we use estimator
$\widehat{FI}(\mathbf{X},\mathcal{D})$ in (\ref{estim1}). In Table
\ref{tab1} are the obtained estimates, as well as, the estimates of
the extremal coefficient ($\widehat{\epsilon}^G$) within each group
and in the whole system (we denote the whole system, i.e., the
vector of all observations as ``Global"). The estimates of the FI
within each financial market group and in the whole system are also
presented. One can see that
the USA is the most fragile financial system with
$\widehat{FI}=1.885948905$. Observe also that the FI of the whole
system is almost twice the FI of the system  divided into blocks.

\begin{table}[!hbp]
\begin{center}
\begin{tabular}{rl}
\begin{tabular}{c|c|c|}
& $\widehat{\epsilon}^G$&$\widehat{FI}$ \\
\hline
Europe & 2.243980009& 1.78254707 \\
USA & 1.590711176& 1.885948905 \\
Far East & 1.673156121 & 1.195345715  \\
\hline
Global & 4.017568517 &  2.240160924
\end{tabular}\hspace{-0.5cm}& \hspace{-0.5cm}
\begin{tabular}{|c}
 $\widehat{FI}(\mathbf{X},\mathcal{D})$\\
\hline
   1.370940479\\
\hline
   \\
\\
\\

\end{tabular}
\end{tabular}
\caption{Estimates of the extremal coefficient
($\widehat{\epsilon}^G$) within each block (Europe, USA, Far East)
and in the whole system (Global) comprising the three blocks, as
well as, estimates of the FI ($\widehat{FI}$) within each block, in
the whole system and within the system divided into blocks
($\widehat{FI}(\mathbf{X},\mathcal{D})$). \label{tab1}}
\end{center}
\end{table}

\end{document}